\documentclass[12pt]{article}
\usepackage{verbatim}
\usepackage{fancyhdr}
\usepackage{graphicx}
\usepackage{mathrsfs}%有更多的数学符号（\mathscr P），比如花体
\usepackage{amssymb}
\usepackage{pst-poly}

\newtheorem{cor}{Corollary}[section]
\newtheorem{lem}{Lemma}[section]
\newtheorem{thm}{Theorem}[section]
\newtheorem{cjt}{Conjecture}[section]
\newtheorem{prob}{Problem}[section]

\newenvironment{pf}[1][Proof]{\noindent\textbf{#1.} }{\hfill\rule{1mm}{2mm}}

\makeatletter \@addtoreset{equation}{section} \makeatother

\begin{document}

\title{Roman Bondage Number of a Graph\thanks {The work was supported by NNSF
of China (No. 11071233).}}
\author
{Fu-Tao Hu, Jun-Ming Xu\footnote{Corresponding
author:
xujm@ustc.edu.cn}\\ \\
{\small Department of Mathematics}  \\
{\small University of Science and Technology of China}\\
{\small Hefei, Anhui, 230026, China} }
\date{}
 \maketitle

\setlength{\baselineskip}{24pt}

\begin{quotation}

\textbf{Abstract}: The Roman dominating function on a graph
$G=(V,E)$ is a function $f: V\rightarrow\{0,1,2\}$ such that each
vertex $x$ with $f(x)=0$ is adjacent to at least one vertex $y$ with
$f(y)=2$. The value $f(G)=\sum\limits_{u\in V(G)} f(u)$ is called
the weight of $f$. The Roman domination number $\gamma_{\rm R}(G)$
is defined as the minimum weight of all Roman dominating functions.
This paper defines the Roman bondage number $b_{\rm R}(G)$ of a
nonempty graph $G=(V,E)$ to be the cardinality among all sets of
edges $B\subseteq E$ for which $\gamma_{\rm R}(G-B)>\gamma_{\rm
R}(G)$. Some bounds are obtained for $b_{\rm R}(G)$, and the exact
values are determined for several classes of graphs. Moreover, the
decision problem for $b_{\rm R}(G)$ is proved to be NP-hard even for
bipartite graphs.

\vskip6pt\noindent{\bf Keywords}: Roman domination number, Roman bondage number,NP-hardness.

\noindent{\bf AMS Subject Classification: }\ 05C69

\end{quotation}

\section{Introduction}

In this paper, a graph $G=(V,E)$ is considered as an undirected
graph without loops and multi-edges, where $V=V(G)$ is the vertex
set and $E=E(G)$ is the edge set. For each vertex $x\in V(G)$, let
$N_G(x)=\{y\in V(G): (x,y)\in E(G)\}$, $N_G[x]=N_G(x)\cup \{x\}$.

A subset $S\subseteq V$ is a {\it dominating set} of $G$ if
$N_G[x]\cap S\ne \emptyset$ for every vertex $x$ in $G$. The {\it
domination number} of $G$, denoted by $\gamma(G)$, is the minimum
cardinality of all dominating sets of $G$. The {\it Roman dominating
function} on a graph $G=(V,E)$, proposed by Cockayne {\it et
al.}~\cite{cd04}, is a function $f: V\rightarrow\{0,1,2\}$ such that
each vertex $x$ with $f(x)=0$ is adjacent to at least one vertex $y$
with $f(y)=2$. Let $(V_0, V_1, V_2)$ be the ordered partition of $V$
induced by $f$, where $V_i=\{v\in V|f(v)=i\}$ for $i=0,1,2$. It is
clear that $V_1\cup V_2$ is a dominating set of $G$, called {\it the
Roman dominating set}, denoted by $D_{\rm R}=(V_1,V_2)$. For
$S\subseteq V$, let $f(S)=\sum\limits_{u\in S} f(u)$. The value
$f(V(G))$ is called the {\it weight} of $f$, denoted by $f(G)$. The
{\it Roman domination number}, denoted by $\gamma_{\rm R}(G)$, is
defined as the minimum weight of all Roman dominating functions,
that is,
 $$
 \gamma_{\rm R}(G)=\min \{f(G): f\ {\rm is\ a\ Roman\ dominating\ function\ on}\ G \}.
 $$

It is clear that for a Roman dominating function $f$ on $G$ and a
Roman dominating set $D_{\rm R}$ of $G$, $f(D_{\rm R}) = 2|V_2| +
|V_1|$. If $D_{\rm R}$ is a minimum Roman dominating set of graph
$G$, then $f(D_{\rm R})=\gamma_{\rm R}(G)$. A Roman dominating
function $f$ is called a {\it $\gamma_{\rm R}$-function} if
$f(G)=\gamma_{\rm R}(G)$. It has been showed by Cockayne {\it et
al.}~\cite{cd04} that for any graph $G$, $\gamma(G)\leqslant
\gamma_R(G)\leqslant2\gamma(G)$. A graph $G$ is called to be {\it
Roman} if $\gamma_{\rm R}(G)=2\gamma(G)$. Roman domination numbers
have been studied, for example, in \cite{cd04, ck10, fk09, fy09,
lk08, lkl05, pp02, rs07, sh07, sh10, xc06}.

To measure the vulnerability or the stability of the domination in
an interconnection network under edge failure, Fink et
at.~\cite{fjkr90} proposed the concept of the bondage number in
1990. The {\it bondage number}, denoted by $b(G)$, of $G$ is the
minimum number of edges whose removal from $G$ results in a graph
with larger domination number of $G$.

Analogously, we can define the Roman bondage number. The {\it Roman
bondage number}, denoted by $b_{\rm R}(G)$, of a nonempty graph $G$
is the minimum number of edges whose removal from $G$ results in a
graph with larger Roman domination number. Precisely speaking, the
Roman bondage number
 $$
 b_{\rm R}(G)=\min\{|B|: B\subseteq E(G), \gamma_{\rm R}(G-B)>\gamma_{\rm R}(G)\}.
 $$

An edge set $B$ that $\gamma_{\rm R}(G-B)>\gamma_{\rm R}(G)$ is
called the {\it Roman bondage set} and the minimum one the {\it
minimum Roman bondage set}. In fact, if $B$ is a minimum Roman
bondage set, then $\gamma_{\rm R}(G-B)=\gamma_{\rm R}(G)+1$, because
the removal of one single edge can not increase the Roman domination
number by more than one. If $b_{\rm R}(G)$ does not exist we define
$b_{\rm R}(G)=\infty$.

In this paper, we give an original investigation. Some bounds are
obtained for $b_{\rm R}(G)$, and the exact values are determined for
several classes of graphs. Moreover, the decision problem for
$b_{\rm R}(G)$ is proved to be NP-hard even for bipartite graphs.

In the proofs of our results, when a Roman dominating function of a
graph is constructed, we only give its nonzero value of some
vertices.

\section{Some basic results on $\gamma_{\rm R}$}

For terminology and notation on graph theory not given here, the
reader is referred to \cite{hh197,hh297,x03}.

Let $G=(V,E)$ be a graph and $E_G(x)=\{xy\in E(G): y\in N_G(x)\}$.
For two disjoint nonempty sets $S,T\subset V(G)$, $E_G(S,T)=E(S,T)$
denotes the set of edges between $S$ and $T$. The degree of $x$ is
denoted by $d_G(x)$, which is equal to $|N_G(x)|$, and $n_i$ denotes
the number of vertices of degree $i$ in $G$ for $i=1,2,\cdots
,\Delta(G)$. Denote the maximum and the minimum degree of $G$ by
$\Delta(G)$ and $\delta(G)$, respectively.

The symbols $P_n$ and $C_n$ denote a path and a cycle, respectively,
where $V(P_n)=V(C_n)=\{x_1,x_2,\cdots,x_n\}$, $E(P_n)=\{x_ix_{i+1}:\
i=1,2,\cdots,n\}$ and $E(C_n)=E(P_n)\cup \{x_1x_n\}$.

The Cartesian product graph $G_1\times G_2$ of two graphs
$G_1=(V_1,E_1)$ and $G_2=(V_2,E_2)$ is a graph with vertex-set
$V=V_1\times V_2=\{(x,y) :\ x\in V_1,y\in V_2\}$, and two vertices
$(x_1,y_1)$ and $(x_2,y_2)$ being adjacent if and only if either
$x_1=x_2$, $y_1$ and $y_2$ are adjacent in $G_2$, or $y_1=y_2$,
$x_1$ and $x_2$ are adjacent in $G_1$.

In this section, we recall some basic results on $\gamma_{\rm R}$,
which will be used in our discussion.

\begin{lem}\label{lem2.1}
{\rm (Cockayne et al. \cite{cd04})} For a path $P_n$ and a cycle
$C_n$,
 $$
 \gamma_{\rm R}(P_n)=\gamma_{\rm R}(C_n)=\left\lceil\frac{2n}{3}\right\rceil.
 $$
For a grid graph $P_2\times P_n$,
 $$
 \gamma_{\rm R}(P_2\times P_n)=n+1.
 $$
For a complete $t$-partite graph $K_{m_{1},m_{2},\cdots,m_t}$ with
$1\leq m_1\leq m_2\leq \cdots \leq m_t$ and $t\geq 2$,
 $$
 \gamma_{\rm R}(K_{m_{1},m_{2},\cdots,m_t})=\left\{
 \begin{array}{ll}
 2, & {\rm if}\ m_1=1;\\
 3, & {\rm if}\ m_1=2;\\
 4, & {\rm if}\ m_1\geq 3.
 \end{array}\right.
 $$
\end{lem}

\begin{lem}\label{lem2.2} {\rm (Cockayne et al. \cite{cd04})}\
If $G$ is a graph of order $n$ and contains vertices of degree
$n-1$, then $\gamma_{\rm R}(G)=2$.
\end{lem}

\begin{lem}\label{lem2.3}
Let $G$ be a nonempty graph with order $n\geq 3$, then
$\gamma_{\rm R}(G)=3$ if and only if $\Delta(G)=n-2$.
\end{lem}

\begin{pf}
Assume that $u$ is a vertex of degree $n-2$ and $v$ is the unique
vertex not adjacent to $u$ in $G$. It is easy to verify that
$\gamma_{\rm R}(G)\geq 3$. Let $f$ be a function from $V(G)$ to
$\{0,1,2\}$ subject to
 $$
 f(x)=\left\{\begin{array}{ll}
 2,\ & {\rm if}\ x=u;\\
 1,\ & {\rm if}\ x=v;\\
 0,\ & {\rm otherwise}.
 \end{array}\right.
 $$
Then $f$ is a Roman dominating function of $G$ with $f(G)=3$. Thus,
$\gamma_{\rm R}(G)=3$.

Conversely, assume $\gamma_{\rm R}(G)=3$, then $\Delta(G)\leq n-2$
by Lemma~\ref{lem2.2}. Let $f$ be a $\gamma_{\rm R}$-function of $G$.

If there is no vertex $u$ with $f(u)=2$, then $f(v)=1$ for each
vertex $v\in V(G)$, and so $n=3$ since $f(G)=\gamma_{\rm R}(G)=3$.
Sine $G$ is nonempty and not $K_3$, $G$ consists of $K_2$ and an
isolated vertex. Thus, $\Delta(G)=1=n-2$.

If there is a vertex $u$ with $f(u)=2$, then there is only one
vertex $v\in V(G)$ with $f(v)=1$ since $f(G)=\gamma_{\rm R}(G)=3$.
The other $n-2$ vertices assigned $0$ are all adjacent to $u$. Thus,
$\Delta(G)\geq d_G(u)\geq n-2$ and hence $\Delta(G)=n-2$.
\end{pf}

\begin{lem}\label{lem2.4}
Let $G$ be an {\rm ($n-3$)}-regular graph with order $n\ge 4$. Then
$\gamma_{\rm R}(G)=4$.
\end{lem}

\begin{pf}
Since $G$ is ($n-3$)-regular and $n\ge 4$, $G$ is nonempty. It is
clear that $\gamma_{\rm R}(G)>2$. By Lemma~\ref{lem2.3},
$\gamma_{\rm R}(G)\neq 3$ since $\Delta(G)=n-3$. Then $\gamma_{\rm
R}(G)\geq 4$. For any vertex $x\in V(G)$, let $y,z$ be the only two
vertices not adjacent to $x$ in $G$, let $f(x)=2$ and $f(y)=f(z)=1$.
Then, $f$ is a Roman dominating function of $G$ with $f(G)=4$, hence
$\gamma_{\rm R}(G)\leq 4$. Thus, $\gamma_{\rm R}(G)=4$.
\end{pf}

\begin{lem}\label{lem2.5}
{\rm (Cockayne et al. \cite{cd04})}\ For any graph $G$,
$\gamma(G)\leq \gamma_{\rm R}(G)\leq 2\gamma(G)$.
\end{lem}

\begin{lem}\label{lem2.6}
{\rm (Cockayne et al. \cite{cd04})}\ For any graph $G$ of order n,
$\gamma(G)=\gamma_{\rm R}(G)$ if and only if $G=\bar{K_n}$.
\end{lem}

\begin{lem}\label{lem2.7}
{\rm (Cockayne et al. \cite{cd04})}\ If $G$ is a connected graph of
order n, then $\gamma_{\rm R}(G) =\gamma(G)+1$ if and only if there
is a vertex $v\in V(G)$ of degree $n-\gamma(G)$.
\end{lem}

A graph $G$ is called to be {\it vertex domination-critical} ( {\it
vc-graph} for short) if $\gamma(G-x) < \gamma(G)$ for any vertex $x$
in $G$.

\begin{lem}\label{lem2.8}
{\rm (Brigham, Chinn and Dutton~\cite{bcd88}, 1988)} A graph $G$
with $\gamma(G)=2$ is a vc-graph if and only if $G$ is a complete
graph $K_{2t}\, (t\geq 2)$ with a perfect matching removed.
\end{lem}

A graph $G$ of order $n$ is {\it vertex Roman domination-critical} (
{\it vrc-graph} for short) if $\gamma_{\rm R}(G)\neq n$ and
$\gamma_{\rm R}(G-x) < \gamma_{\rm R}(G)$ for any vertex $x$ in $G$.
For example, for a positive integer $k$, both $C_{3k+1}$ and
$C_{3k+2}$ are vrc-graphs by Lemma~\ref{lem2.1}. From the
definition, it is clear that $\gamma_{\rm R}(G)\geq 3$ if $G$ is a
vrc-graph with order at least 3.

\begin{lem}\label{lem2.9}
If $G$ is a vrc-graph with $\gamma_{\rm R}(G)=3$, then $G$ is a
vc-graph with $\gamma(G)=2$.
\end{lem}

\begin{pf}
Let $G$ be a vrc-graph with $\gamma_{\rm R}(G)=3$. From the
definition of vrc-graph, $|V(G)|>\gamma_{\rm R}(G)=3$. By
Lemma~\ref{lem2.3}, $\Delta(G)=|V(G)|-2$ and hence $\gamma(G)=2$.
For any vertex $x$, if $\gamma_{\rm R}(G-x) < \gamma_{\rm R}(G)=3$,
then $\gamma_{\rm R}(G-x)=2$ since $G-x$ is nonempty. By
Lemma~\ref{lem2.2}, $G-x$ contains vertices of degree $|V(G-x)|-1$
and, hence, $\gamma(G-x)=1$, which implies that $G$ is a vc-graph.
\end{pf}

\section{The exact values of $b_{\rm R}$ for some graphs}

\begin{lem}\label{lem3.1}
Let $G$ be a graph with order $n\geq 3$ and
%$\gamma_{\rm R}(G)=2$,
$t$ be the number of vertices of degree $n-1$ in $G$. If $t\geq 1$
then $b_{\rm R}(G)=\lceil\frac{t}{2}\rceil$.
\end{lem}

\begin{pf}
Let $H$ be a spanning subgraph of $G$ obtained by removing fewer
than $\lceil\frac{t}{2}\rceil$ edges from $G$. Then $H$ contains
vertices of degree $n-1$ and, hence, $\gamma_{\rm
R}(H)=2=\gamma_{\rm R}(G)$ by Lemma~\ref{lem2.2}, which implies
$b_{\rm R}(G)\geq \lceil\frac{t}{2}\rceil$.

Since $G$ contains $t$ vertices of degree $n-1$, it contains a
complete subgraph $K_t$ induced by these $t$ vertices. We can remove
$\lceil\frac{t}{2}\rceil$ edges such that no vertices have degree
$n-1$ and, hence, $\gamma_{\rm R}(H)\geq 3>2=\gamma_{\rm R}(G)$
since $n\geq 3$. Thus $b_{\rm R}(G)\leq \lceil\frac{t}{2}\rceil$,
whence $b_{\rm R}(G)=\lceil\frac{t}{2}\rceil$.
\end{pf}

\begin{cor}\label{cor3.1}
For a complete graph $K_n$ ($n\geq 3$), $b_{\rm
R}(K_n)=\lceil\frac{n}{2}\rceil$.
\end{cor}

\begin{thm}\label{thm3.2} For a path $P_n$ with $n\geq 3$,
 $$
b_{\rm R}(P_n)= \left\{ \begin{array}{ll}
 1, & {\rm if}\ n\equiv 0,1\, ({\rm mod}\, 3);\\
 2, & {\rm otherwise}.
 \end{array} \right.
 $$
\end{thm}
\begin{pf}
Let $P_n=(x_1,x_2,\ldots,x_n)$ be a path. By Lemma~\ref{lem2.1},
$\gamma_{\rm R}(P_n)=\lceil\frac{2n}{3}\rceil$.

If $n\equiv0,1\, ({\rm mod}\, 3)$, then
 $$
 \gamma_{\rm R}(P_n-x_2x_3)=2+\left\lceil\frac{2(n-2)}{3}\right\rceil
 =1+\left\lceil\frac{2n-1}{3}\right\rceil=1+\gamma_{\rm R}(P_n),
$$
and hence $b_{\rm R}(P_n)\leq 1$, whence $b_{\rm
R}(P_n)=1$.

If $n\equiv 2\, ({\rm mod}\, 3)$, then for any edge $e=x_ix_{i+1}\in
E(P_n)$,
 $$
 \gamma_{\rm R}(P_n-e)=\left\lceil\frac{2i}{3}\right\rceil+\left\lceil\frac{2(n-i)}{3}\right\rceil \leq
 \left\lceil\frac{2(n-i)+2i+2}{3}\right\rceil=\left\lceil\frac{2n}{3}\right\rceil=\gamma_{\rm R}(P_n),
 $$
and hence $b_{\rm R}(P_n)\geq 2$. Since
 $$
 \gamma_{\rm R}(P_n-x_2x_3-x_4x_5)=2+2+\left\lceil\frac{2(n-4)}{3}\right\rceil
 =1+\left\lceil\frac{2n+1}{3}\right\rceil\geq 1+\gamma_{\rm R}(P_n),
 $$
we have $b_{\rm R}(P_n)\leq 2$, whence $b_{\rm R}(P_n)=2$.
\end{pf}

\begin{cor}\label{cor3.2} For a cycle $C_n$ with $n\geq 3$,
 $$
b_{\rm R}(C_n)= \left\{ \begin{array}{ll}
 2, & {\rm if}\ n=0,1\, ({\rm mod}\, 3);\\
 3, & {\rm otherwise}.
 \end{array} \right.
 $$
\end{cor}

\begin{lem}\label{lem3.2}
Let $P_n=(x_1,x_2,\ldots,x_n)$ be a path, and use $u_{i,j}$ to
denote the vertex $(x_i,x_j)$ in $P_2\times P_n$, where $1\le i\le
2$ and $1\le j\le n$. Then there exists a $\gamma_{\rm R}$-function
$f$ on $P_2\times P_n$ such that $f(u_{1,1})=2$ or $f(u_{2,1})=2$ or
$f(u_{1,n})=2$ or $f(u_{2,n})=2$.
\end{lem}

\begin{pf}
Without loss of generality, we only need to find a $\gamma_{\rm
R}$-function $f$ on $P_2\times P_n$ with $f(u_{1,1})=2$. Define a
Roman dominating function $f$ as follows. For each non-negative
integer $i$ with $1+4i\leq n$, let $f(u_{1,1+4i})=2$, and for each
non-negative integer $j$ with $3+4j\leq n$, let $f(u_{2,3+4j})=2$.
If $n\equiv 0\,({\rm mod}\,4)$, let $f(u_{1,n})=1$, and if $n\equiv
2\,({\rm mod}\,4)$, let $f(u_{2,n})=1$. Then $f(P_2\times P_n)=n+1$
and, hence by Lemma~\ref{lem2.1}, $f$ is a $\gamma_{\rm R}$-function
with $f(u_{1,1})=2$.
\end{pf}

\begin{thm}\label{thm3.2}
$b_{\rm R}(P_2\times P_n)=2$ for $n\geq 2$.
\end{thm}

\begin{pf}
By Lemma~\ref{lem2.1}, we have $\gamma_{\rm R}(P_2\times P_n)=n+1$.
Since $\gamma_{\rm R}(P_2\times
P_n-u_{1,1}u_{1,2}-u_{2,1}u_{2,2})=2+\gamma_{\rm R}(P_2\times
P_{n-1})=n+2$, we have $b_{\rm R}(P_2\times P_n)\leq 2$. Next we
prove that $\gamma_{\rm R}(P_2\times P_n-e)\leq \gamma_{\rm
R}(P_2\times P_n)$ for any edge $e\in E(P_2\times P_n)$.

Suppose that $e$ is incident with some vertex in
$\{u_{1,1},u_{2,1},u_{1,n},u_{2,n}\}$. Without loss of generality
let $e$ be incident with $u_{1,1}$. By Lemma~\ref{lem3.2}, there
exists a $\gamma_{\rm R}$-function $f$ on $P_2\times (P_n-P_1)$ such
that $f(u_{2,2})=2$. Denote $f(u_{1,1})=1$ and then $f$ is a Roman
dominating function of $P_2\times P_n-e$ with $f(P_2\times P_n-e)=n+1$, thus
$\gamma_{\rm R}(P_2\times P_n-e)\leq \gamma_{\rm R}(P_2\times P_n)$.

Suppose that $e$ is incident with some vertex in $\{u_{i,j}:\ 1\le
i\le 2, 2\leq j \leq n-1\}\setminus\{u_{1,1},u_{2,1},
u_{1,n},u_{2,n}\}$. Without loss of generality let $e$ be incident
with $u_{1,j}$ and not incident with $u_{1,j-1}$. By
Lemma~\ref{lem3.2}, there exists a $\gamma_{\rm R}$-function $f_1$
on $P_2\times P_{j-1}$ with $f_1(u_{1,j-1})=2$ and a $\gamma_{\rm
R}$-function $f_2$ on $P_2\times (P_n-P_j)$ with $f_2(u_{2,j+1})=2$.
Then $f=f_1\cup f_2$ is a Roman dominating function on $P_2\times
P_n-e$ with $f(P_2\times P_n-e)=n+1$, thus $\gamma_{\rm R}(P_2\times
P_n-e)\leq \gamma_{\rm R}(P_2\times P_n)$.

The above two cases yield that $b_{\rm R}(P_2\times P_n)\geq 2$ and,
hence, $b_{\rm R}(P_2\times P_n)=2$. The lemma follows.
\end{pf}

\section{Complexity of Roman bondage number}

In this section, we will show that the Roman bondage number problem
is NP-hard and the Roman domination number problem is NP-complete
even for bipartite graphs. We first state the problem as the
following decision problem.

\begin{center}
\begin{minipage}{130mm}
\setlength{\baselineskip}{24pt}

\vskip6pt\noindent {\bf Roman bondage number problem (RBN):}

\noindent {\bf Instance:}\ {\it A nonempty bipartite graph $G$ and a positive
integer $k$.}

\noindent {\bf Question:}\ {\it Is $b_{\rm R}(G)\le k$?}

\end{minipage}
\end{center}

\begin{center}
\begin{minipage}{130mm}
\setlength{\baselineskip}{24pt}

\vskip6pt\noindent {\bf Roman domination number problem (RDN):}

\noindent {\bf Instance:}\ {\it A nonempty bipartite graph $G$ and a positive
integer $k$.}

\noindent {\bf Question:}\ {\it Is $\gamma_{\rm R}(G)\le k$?}

\end{minipage}
\end{center}

Following Garey and Johnson's techniques for proving NP-completeness
given in~\cite{gj79}, we prove our results by describing a
polynomial transformation from the known-well NP-complete problem:
3SAT. To state 3SAT, we recall some terms.

Let $U$ be a set of Boolean variables. A {\it truth assignment} for
$U$ is a mapping $t: U\to\{T,F\}$. If $t(u)=T$, then $u$ is said to
be ``\,true" under $t$; If $t(u)=F$, then $u$ is said to be
``\,false" under $t$. If $u$ is a variable in $U$, then $u$ and
$\bar{u}$ are {\it literals} over $U$. The literal $u$ is true under
$t$ if and only if the variable $u$ is true under $t$; the literal
$\bar{u}$ is true if and only if the variable $u$ is false.

A {\it clause} over $U$ is a set of literals over $U$. It represents
the disjunction of these literals and is {\it satisfied} by a truth
assignment if and only if at least one of its members is true under
that assignment. A collection $\mathscr C$ of clauses over $U$ is
{\it satisfiable} if and only if there exists some truth assignment
for $U$ that simultaneously satisfies all the clauses in $\mathscr
C$. Such a truth assignment is called a {\it satisfying truth
assignment} for $\mathscr C$. The 3SAT is specified as follows.

\begin{center}
\begin{minipage}{130mm}
\setlength{\baselineskip}{24pt}

\vskip6pt\noindent {\bf $3$-satisfiability problem (3SAT):}

\noindent {\bf Instance:}\ {\it A collection
$\mathscr{C}=\{C_1,C_2,\ldots,C_m\}$ of clauses over a finite set
$U$ of variables such that $|C_j| =3$ for $j=1, 2,\ldots,m$.}

\noindent {\bf Question:}\ {\it Is there a truth assignment for $U$
that satisfies all the clauses in $\mathscr{C}$?}

\end{minipage}
\end{center}

\begin{thm} \textnormal{(Theorem 3.1 in~\cite{gj79})}
3SAT is NP-complete.
\end{thm}

\begin{thm}\label{thm4.1}
RBN is NP-hard even for bipartite graphs.
\end{thm}

\begin{pf}
The transformation is from 3SAT. Let $U=\{u_1,u_2,\ldots,u_n\}$ and
$\mathscr{C}=\{C_1,C_2,\ldots,$ $C_m\}$ be an arbitrary instance of
3SAT. We will construct a bipartite graph $G$ and choose an integer
$k$ such that $\mathscr{C}$ is satisfiable if and only if $b_{\rm
R}(G)\leq k$. We construct such a graph $G$ as follows.

For each $i=1,2,\ldots,n$, corresponding to the variable $u_i\in U$,
associate a graph $H_i$ with vertex set
$V(H_i)=\{u_i,\bar{u}_i,v_i,v_i',x_i,y_i,z_i,w_i\}$ and edge set
$E(H_i)=\{u_iv_i,u_iz_i,\bar{u}_iv_i',\\
\bar{u}_iz_i,y_iv_i,y_iv_i',y_iz_i,w_iv_i,w_iv_i',w_iz_i,
x_iv_i,x_iv_i'\}$. For each $j=1,2,\ldots,m$, corresponding to the
clause $C_j=\{p_j,q_j,r_j\}\in \mathscr{C}$, associate a single
vertex $c_j$ and add edge set $E_j=\{c_jp_j, c_jq_j,c_jr_j\}$, $1\le
j\le m$. Finally, add a path $P=s_1s_2s_3$, join $s_1$ and $s_3$ to
each vertex $c_j$ with $1\le j\le m$ and set $k=1$.

Figure~\ref{f1} shows an example of the graph obtained when
$U=\{u_1,u_2,u_3,u_4\}$ and $\mathscr{C}=\{C_1,C_2,C_3\}$, where
$C_1=\{u_1,u_2,\bar{u}_3\}, C_2=\{\bar{u}_1,u_2,u_4\},
C_3=\{\bar{u}_2,u_3,u_4\}$.

\begin{figure}[h]
\begin{center}
\begin{pspicture}(-5,-1.1)(5,7.7)

\cnode*(0,-.6){3pt}{s2}\rput(0,-1){$s_2$}
\cnode(-1,0){3pt}{s1}\rput(-1,-.4){$s_1$}
\cnode(1,0){3pt}{s3}\rput(1.1,-.4){$s_3$} \ncline{s2}{s1}
\ncline{s2}{s3}

\cnode(0,1.9){3pt}{c2}\rput(0,2.3){$c_2$}
\ncline{c2}{s1} \ncline{c2}{s3}
\cnode(-2.5,2){3pt}{c1}\rput(-2.8,1.7){$c_1$}
\ncline{c1}{s1} \ncline{c1}{s3}
\cnode(2.5,2){3pt}{c3}\rput(2.8,1.7){$c_3$}
\ncline{c3}{s1} \ncline{c3}{s3}

\cnode(-4.5,4){3pt}{u1}\rput(-4.7,4.3){$u_1$}
\cnode(-3,4){3pt}{u1'}\rput(-2.8,4.3){$\bar{u}_1$}

\cnode*(-2,4){3pt}{u2}\rput(-2.2,4.3){$u_2$}
\cnode(-0.5,4){3pt}{u2'}\rput(-0.3,4.3){$\bar{u}_2$}

\cnode(0.5,4){3pt}{u3}\rput(0.3,4.3){$u_3$}
\cnode*(2,4){3pt}{u3'}\rput(2.2,4.3){$\bar{u}_3$}

\cnode*(3,4){3pt}{u4}\rput(2.8,4.3){$u_4$}
\cnode(4.5,4){3pt}{u4'}\rput(4.7,4.3){$\bar{u}_4$}

\cnode(-4.5,6){3pt}{v1}\rput(-4.7,6.3){$v_1$}
\cnode(-3,6){3pt}{v1'}\rput(-2.8,6.3){$v_1'$}
\ncline{u1}{v1} \ncline{u1'}{v1'}
\cnode(-2,6){3pt}{v2}\rput(-2.2,6.3){$v_2$}
\cnode*(-0.5,6){3pt}{v2'}\rput(-0.3,6.3){$v_2'$}
\ncline{u2}{v2} \ncline{u2'}{v2'}
\cnode*(0.5,6){3pt}{v3}\rput(0.3,6.3){$v_3$}
\cnode(2,6){3pt}{v3'}\rput(2.2,6.3){$v_3'$}
\ncline{u3}{v3} \ncline{u3'}{v3'}
\cnode(3,6){3pt}{v4}\rput(2.8,6.3){$v_4$}
\cnode*(4.5,6){3pt}{v4'}\rput(4.7,6.3){$v_4'$}
\ncline{u4}{v4} \ncline{u4'}{v4'}

\cnode(-3.75,5){3pt}{w1}\rput(-3.45,4.7){$w_1$}
\ncline{w1}{v1} \ncline{w1}{v1'}
\cnode(-1.25,5){3pt}{w2}\rput(-.95,4.7){$w_2$}
\ncline{w2}{v2} \ncline{w2}{v2'}
\cnode(1.25,5){3pt}{w3}\rput(1.55,4.7){$w_3$}
\ncline{w3}{v3} \ncline{w3}{v3'}
\cnode(3.75,5){3pt}{w4}\rput(4.05,4.7){$w_4$}
\ncline{w4}{v4} \ncline{w4}{v4'}

\cnode*(-3.75,4){3pt}{z1}\rput(-3.75,3.7){$z_1$}
\ncline{z1}{u1} \ncline{z1}{u1'} \ncline{w1}{z1}
\cnode(-1.25,4){3pt}{z2}\rput(-1.25,3.7){$z_2$}
\ncline{z2}{u2} \ncline{z2}{u2'}  \ncline{w2}{z2}
\cnode(1.25,4){3pt}{z3}\rput(1.45,4.3){$z_3$}
\ncline{z3}{u3} \ncline{z3}{u3'}  \ncline{w3}{z3}
\cnode(3.75,4){3pt}{z4}\rput(3.75,3.7){$z_4$}
\ncline{z4}{u4} \ncline{z4}{u4'}  \ncline{w4}{z4}

\cnode(-3.75,6){3pt}{y1}\rput(-3.75,6.3){$y_1$}
\ncline{y1}{v1} \ncline{y1}{v1'}  \nccurve[angleA=-135,angleB=135]{y1}{z1}
\cnode(-1.25,6){3pt}{y2}\rput(-1.25,6.3){$y_2$}
\ncline{y2}{v2} \ncline{y2}{v2'}  \nccurve[angleA=-135,angleB=135]{y2}{z2}
\cnode(1.25,6){3pt}{y3}\rput(1.25,6.3){$y_3$}
\ncline{y3}{v3} \ncline{y3}{v3'}  \nccurve[angleA=-135,angleB=135]{y3}{z3}
\cnode(3.75,6){3pt}{y4}\rput(3.75,6.3){$y_4$}
\ncline{y4}{v4} \ncline{y4}{v4'}  \nccurve[angleA=-135,angleB=135]{y4}{z4}

\cnode*(-3.75,7){3pt}{x1}\rput(-3.75,7.3){$x_1$}
\ncline{x1}{v1} \ncline{x1}{v1'}
\cnode(-1.25,7){3pt}{x2}\rput(-1.25,7.3){$x_2$}
\ncline{x2}{v2} \ncline{x2}{v2'}
\cnode(1.25,7){3pt}{x3}\rput(1.25,7.3){$x_3$}
\ncline{x3}{v3} \ncline{x3}{v3'}
\cnode(3.75,7){3pt}{x4}\rput(3.75,7.3){$x_4$}
\ncline{x4}{v4} \ncline{x4}{v4'}

\ncline{c1}{u1} \ncline{c1}{u2} \ncline{c1}{u3'}
\ncline{c2}{u1'} \ncline{c2}{u2} \ncline{c2}{u4}
\ncline{c3}{u2'} \ncline{c3}{u3} \ncline{c3}{u4}
\end{pspicture}
\caption{\label{f1}\footnotesize An instance of the Roman bondage
number problem resulting from an instance of 3SAT. Here $k=1$ and
$\gamma_{\rm R}(G)=18$, where the bold vertex $w$ means a Roman
dominating function with $f(w)=2$.}
\end{center}
\end{figure}
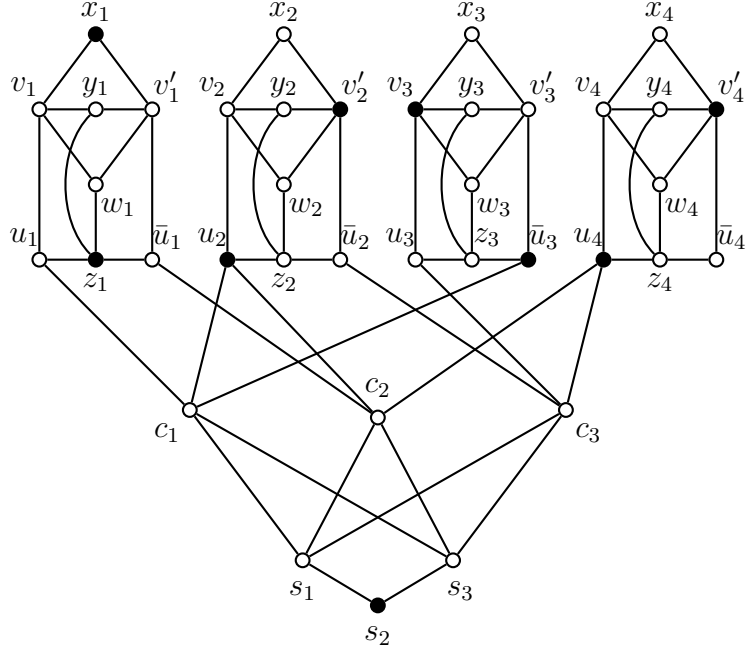

To prove that this is indeed a transformation, we only need to show
that $b_{\rm R}(G)=1$ if and only if there is a truth assignment for
$U$ that satisfies all clauses in $\mathscr{C}$. This aim can be
obtained by proving the following four claims.

\begin{description}

\item [Claim 4.1]
{\it $\gamma_{\rm R}(G)\geq 4n+2$. Moreover, if $\gamma_{\rm
R}(G)=4n+2$, then for any $\gamma_{\rm R}$-function $f$ on $G$,
$f(H_i)=4$ and at most one of $f(u_i)$ and $f(\bar{u}_i)$ is 2 for
each $i$, $f(c_j)=0$ for each $j$ and $f(s_2)=2$.}

\begin{pf}
Let $f$ be a $\gamma_{\rm R}$-function of $G$, and let
$H_i'=H_i-u_i-\bar{u}_i$.

If $f(u_i)=2$ and $f(\bar{u}_i)=2$, then $f(H_i)\geq 4$. Assume
either $f(u_i)=2$ or $f(\bar{u}_i)=2$, if $f(x_i)=0$ or $f(y_i)=0$,
then there is at least one vertex $t$ in $\{v_i, \bar{v}_i,z_i\}$
such that $f(t)=2$. And hence $f(H_i')\ge 2$. Thus, $f(H_i)\geq 4$.

If $f(u_i)\neq 2$ and $f(\bar{u}_i)\neq 2$, let $f'$ be a
restriction of $f$ on $H_i'$, then $f'$ is a Roman dominating
function of $H_i'$, and $f'(H_i')\geq \gamma_{\rm R}(H_i')$. Since
the maximum degree of $H_i'$ is $V(H_i')-3$, by Lemma~\ref{lem2.3},
$\gamma_{\rm R}(H_i')>3$ and hence $f'(H_i')\geq 4$ and $f(H_i)\geq
4$. If $f(s_1)=0$ or $f(s_3)=0$, then there is at least one vertex
$t$ in $\{c_1, \cdots,c_m, s_2\}$ such that $f(t)=2$. Then
$f(N_G[V(P)])\geq 2$, and hence $\gamma_{\rm R}(G)\geq 4n+2$.

Suppose that $\gamma_{\rm R}(G)=4n+2$, then $f(H_i)=4$ and since
$f(N_G[x_i])\ge 1$, at most one of $f(u_i)$ and $f(\bar{u}_i)$ is 2
for each $i=1,2,\ldots,n$, while $f(N_G[V(P)])=2$. Then we have
$f(s_2)=2$ since $f(N_G[s_2])\ge 1$. Consequently, $f(c_j)=0$ for
each $j=1,2,\ldots,m$.
\end{pf}

\item [Claim 4.2]
{\it $\gamma_{\rm R}(G)=4n+2$ if and only if $\mathscr{C}$ is satisfiable.}

\begin{pf}
Suppose that $\gamma_{\rm R}(G)=4n+2$ and let $f$ be a $\gamma_{\rm
R}$-function of $G$. By Claim 4.1, at most one of $f(u_i)$ and
$f(\bar{u}_i)$ is 2 for each $i=1,2,\ldots,n$. Define a mapping $t:
U\to \{T,F\}$ by
\begin{equation}\label{e4.1}
 t(u_i)=\left\{
 \begin{array}{l}
 T \ \ {\rm if}\ f(u_i)=2\ {\rm or}\ f(u_i)\neq 2\ {\rm and} f(\bar{u}_i)\neq 2, \\
 F \ \ {\rm if}\ f(\bar{u}_i)=2.
\end{array}
 \right.
 \ i=1,2,\ldots,n.
 \end{equation}

We now show that $t$ is a satisfying truth assignment for
$\mathscr{C}$. It is sufficient to show that every clause in
$\mathscr{C}$ is satisfied by $t$. To this end, we arbitrarily
choose a clause $C_j\in\mathscr{C}$ with $1\le j\le m$.

By Claim 4.1, $f(c_j)=f(s_1)=f(s_3)=0$. There exists some $i$ with
$1\le i\le n$ such that $f(u_i)=2$ or $f(\bar{u}_i)=2$ where $c_j$
is adjacent to $u_i$ or $\bar{u}_i$. Suppose that $c_j$ is adjacent
to $u_i$ where $f(u_i)=2$. Since $u_i$ is adjacent to $c_j$ in $G$,
the literal $u_i$ is in the clause $C_j$ by the construction of $G$.
Since $f(u_i)=2$, it follows that $t(u_i)=T$ by (\ref{e4.1}), which
implies that the clause $C_j$ is satisfied by $t$. Suppose that
$c_j$ is adjacent to $\bar{u}_i$ where $f(\bar{u}_i)=2$. Since
$\bar{u}_i$ is adjacent to $c_j$ in $G$, the literal $\bar{u}_i$ is
in the clause $C_j$. Since $f(\bar{u}_i)=2$, it follows that
$t(u_i)=F$ by (\ref{e4.1}). Thus, $t$ assigns $\bar{u}_i$ the truth
value $T$, that is, $t$ satisfies the clause $C_j$. By the
arbitrariness of $j$ with $1\le j\le m$, we show that $t$ satisfies
all the clauses in $\mathscr{C}$, that is, $\mathscr{C}$ is
satisfiable.

Conversely, suppose that $\mathscr{C}$ is satisfiable, and let $t:
U\to \{T,F\}$ be a satisfying truth assignment for $\mathscr{C}$.
Create a function $f$ on $V(G)$ as follows: if $t(u_i)=T$, then let
$f(u_i)=f(v_i')=2$, and if $t(u_i)=F$, then let
$f(\bar{u}_i)=f(v_i)=2$. Let $f(s_2)=2$. Clearly, $f(G)=4n+2$. Since
$t$ is a satisfying truth assignment for $\mathscr{C}$, for each
$j=1,2,\ldots,m$, at least one of literals in $C_j$ is true under
the assignment $t$. It follows that the corresponding vertex $c_j$
in $G$ is adjacent to at least one vertex $w$ with $f(w)=2$ since
$c_j$ is adjacent to each literal in $C_j$ by the construction of
$G$. Thus $f$ is a Roman dominating function of $G$, and so
$\gamma_{\rm R}(G)\leq f(G)= 4n+2$. By Claim 4.1, $\gamma_{\rm
R}(G)\geq 4n+2$, and so $\gamma_{\rm R}(G)=4n+2$.
\end{pf}

\item [Claim 4.3]
{\it $\gamma_{\rm R}(G-e)\leq 4n+3$ for any $e\in E(G)$.}

\begin{pf}
For any edge $e\in E(G)$, it is sufficient to construct a Roman
dominating function $f$ with weight $4n+3$ of $G$. We first assume
$e\in E_G(s_1)$ or $e\in E_G(s_3)$ or $e\in E_G(c_j)$ for each
$j=1,2,\ldots,m$, without loss of generality let $e\in E_G(s_1)$ or
$e=c_ju_i$ or $e=c_j\bar{u}_i$. Let $f(s_3)=2, f(s_1)=1$ and
$f(u_i)=f(v_i')=2$ for each $i=1,2,\ldots,n$. For the edge $e\notin
E_G(u_i)$ and $e\notin E_G(v_i')$ or $e=\bar{u}_iz_i$, let
$f(s_1)=2, f(s_3)=1$ and $f(u_i)=f(v_i')=2$. For the edge $e\notin
E(\bar{u}_i)$ and $e\notin E(v_i)$ or $e=u_iz_i$, let $f(s_1)=2,
f(s_3)=1$ and $f(\bar{u}_i)=f(v_i)=2$. If $e=u_iv_i$ or
$e=\bar{u}_iv_i'$, let $f(s_1)=2, f(s_3)=1$ and $f(x_i)=f(z_i)=2$.
Then $f$ is a Roman dominating function of $G-e$ with $f(G-e)=4n+3$
and hence $\gamma_{\rm R}(G-e)\leq 4n+3$.
\end{pf}

\item [Claim 4.4]
{\it $\gamma_{\rm R}(G)=4n+2$ if and only if $b_{\rm R}(G)=1$.}

\begin{pf}
Assume $\gamma_{\rm R}(G)=4n+2$ and consider the edge $e=s_1s_2$.
Suppose $\gamma_{\rm R}(G)=\gamma_{\rm R}(G-e)$. Let $f'$ be a
$\gamma_{\rm R}$-function of $G-e$. It is clear that $f'$ is also a
$\gamma_{\rm R}$-function on $G$. By Claim 4.1 we have  $f'(c_j)=0$
for each $j=1,2,\ldots,m$ and $f'(s_2)=2$. But then
$f'(N_{G-e}[s_1])=0$, a contradiction. Hence, $\gamma_{\rm
R}(G)<\gamma_{\rm R}(G-e)$, and so $b_{\rm R}(G)=1$.

Now, assume $b_{\rm R}(G)=1$. By Claim 4.1, we have that
$\gamma_{\rm R}(G)\geq 4n+2$. Let $e'$ be an edge such that
$\gamma_{\rm R}(G)<\gamma_{\rm R}(G-e')$. By Claim 4.3, we have that
$\gamma_{\rm R}(G-e')\leq 4n+3$. Thus, $4n + 2\leq \gamma_{\rm
R}(G)< \gamma_{\rm R}(G-e')\leq 4n+3$, which yields $\gamma_{\rm
R}(G)=4n+2$.
\end{pf}
\end{description}

By Claim 4.2 and Claim 4.4, we prove that $b_{\rm R}(G)=1$ if and
only if there is a truth assignment for $U$ that satisfies all
clauses in $\mathscr{C}$. Since the construction of the Roman
bondage number instance is straightforward from a $3$-satisfiability
instance, the size of the Roman bondage number instance is bounded
above by a polynomial function of the size of $3$-satisfiability
instance. It follows that this is a polynomial reduction.

The theorem follows.
\end{pf}

\begin{cor}
Roman domination number problem is NP-complete even for bipartite
graphs.
\end{cor}

\begin{pf}
It is easy to see that the Roman bondage problem is in NP since a
nondeterministic algorithm need only guess a vertex set pair
$(V_1,V_2)$ with $|V_1|+2|V_2|\le k$ and check in polynomial time
whether that for any vertex $u\in V\setminus(V_1\cup V_2)$ whether
there is a vertex in $V_2$ adjacent to $u$ for a given nonempty
graph $G$.

We use the same method as Theorem~\ref{thm4.1} to prove this
conclusion. We construct the same graph $G$ but does not contain the
path $P$. We set $k=4n$, then use the same methods as Claim 4.1 and
4.2, we have that $\gamma_{\rm R}(G)=4n$ if and only if
$\mathscr{C}$ is satisfiable.
\end{pf}

\section{General bounds}

\begin{lem}\label{lem5.1}
Let $H$ be a spanning subgraph obtained by removing
$k$ edges from a graph $G$. Then $b_{\rm R}(G)\leq b_{\rm R}(H)+k$.
\end{lem}
\begin{pf}
Let $B= E(G)\setminus E(H)$ and $B'$ be a minimum Roman bondage set
of $H$. Then $|B|=k$, $|B'|=b_{\rm R}(H)$ and $\gamma_{\rm
R}(H-B')>\gamma_{\rm R}(H)$. Let $f: V\rightarrow\{0,1,2\}$ be a
Roman dominating function on $H$ with $f(H)=\gamma_{\rm R}(H)$. Then
each vertex $x$ with $f(x)=0$ is adjacent to at least one vertex $y$
with $f(y)=2$ in $H$, and so is in $G$ since $H=G-B$, which implies
that $f$ is a Roman dominating function of $G$, and so $f(G)\geq
\gamma_{\rm R}(G)$. It follows that $\gamma_{\rm
R}(G-B-B')=\gamma_{\rm R}(H-B')>\gamma_{\rm R}(H)\geq \gamma_{\rm
R}(G)$ and, hence, $b_{\rm R}(G)\leq |B|+|B'|=b_{\rm R}(H)+k$.
\end{pf}

\begin{thm}\label{thm5.1}
$b_{\rm R}(G)\leq d_G(x)+d_G(y)+d_G(z)-|N_G(y)\cap N_G(\{x,z\})|-3$
for any path $(x,y,z)$ of length 2 in a graph $G$.
\end{thm}

\begin{pf}
Let $F_y=\{(y,u)\in E(G): u\in N_G(y)\cap N_G(\{x,z\})\}$,
$B=E_G(x)\cup E_G(z)\cup (E_G(y)\setminus F_y)$. Then
 $$
|B|=d_G(x)+d_G(y)+d_G(z)-|N_G(y)\cap N_G(\{x,z\})|-2.
 $$
Let $H=G-B+yz$. Then $x$ is an isolated vertex and $z$ is a vertex
of degree 1 which is only adjacent to $y$ in $H$. Let $f$ be a
minimum Roman dominating function of $H$, then $f(x)=1$ and $1\leq
f(y)+f(z)\leq 2$.

If $f(y)+f(z)=2$, then let $f'=f$ except $f'(x)=0$, $f'(y)=2$ and
$f'(z)=0$. Clearly, $f'$ is a Roman dominating function of $G$ with
$f'(G)<f(H)$ and, hence, $b_{\rm R}(G)\leq |B|-1$.

If $f(y)+f(z)=1$, then $f(y)=0$ and $f(z)=1$. There is an edge
$(u,y)\in F_y$ with $f(u)=2$. Let $f'=f$ except $f'(x)=0$ if $u\in
N_G(x)$ or $f'(z)=0$ if $u\in N_G(z)\setminus N_G(x)$. Then $f'$ is
a Roman dominating function of $G$ with $f'(G)<f(H)$, and hence
$b_{\rm R}(G)\leq |B|-1$.
\end{pf}

\begin{thm}\label{thm5.2}
$b_{\rm R}(G)\leq d_G(x)+d_G(y)+d_G(z)-|N_G(y)\cap
N_G(\{x,z\})|-|N_G(x)\cap N_G(z)|-1$ for any path $(x,y,z)$ of length 2 in a graph $G$.
\end{thm}

\begin{pf}
Let $F_y=\{(y,u)\in E(G): u\in N_G(y)\cap N_G(\{x,z\})\}$ and
$F_z=\{(z,u)\in E(G): u\in (N_G(z)\cap N_G(x))\}$, $B=E_G(x)\cup
(E_G(z)\setminus F_z)\cup (E_G(y)\setminus F_y)$ and $H=G-B$. Then
$x$ is an isolated vertex in $H$. Let $f$ be a minimum Roman
dominating function of $H$, then $f(x)=1$. We will construct a Roman
dominating function $f'$ of $G$ with $f'(G)<f(H)$.

If $f(z)=0$, then there is an edge $(z,s)\in F_z$ with $f(s)=2$.
Thus, if $f(y)=2$ or $f(z)=0$, let $f'=f$ except $f'(x)=0$. In the
following, let $f(y)\neq 2$ and $f(z)\neq 0$.

If $f(y)=0$. Then there is a vertex $s\in F_y$
such that $f(s)=2$. If $s\in N_G(x)$, let $f'=f$ except
$f'(x)=0$. If $s\in N_G(z)\setminus N_G(x)$, let
$f'=f$ except $f'(z)=0$.

If $f(y)=1$. If $f(z)=1$, let $f'=f$ except $f'(x)=f'(z)=0$ and
$f'(y)=2$. If $f(z)=2$, let $f'=f$ except $f'(y)=0$.

Then $f'$ is a Roman dominating function of $G$ with $f'(G)<f(H)$,
and hence $b_{\rm R}(G)\leq |B|\le d_G(x)+d_G(y)+d_G(z)-|N_G(y)\cap
N_G(\{x,z\})|-|N_G(x)\cap N_G(z)|-1$.
\end{pf}

\begin{cor}\label{cor5.1}
$b_{\rm R}(G)\leq \min\{d_G(x)+d_G(y)+d_G(z)-|N_G(y)\cap
N_G(\{x,z\})|-3, d_G(x)+d_G(y)+d_G(z)-|N_G(y)\cap
N_G(\{x,z\})|-|N_G(x)\cap N_G(z)|-1\}$ for any path $(x,y,z)$ of length 2 in a graph $G$.
\end{cor}

\begin{cor}\label{cor5.2}
$b_{\rm R}(G)\leq 2\Delta(G)+\delta(G)-3$ for any graph with
diameter at least two.
\end{cor}

\begin{cor}\label{cor5.3}
For any tree $T$ of order at least 3, then $b_{\rm R}(T)\leq
\Delta(T)$.
\end{cor}

\begin{pf}
If there is a vertex $x$ adjacent to at least two vertices of degree
one in $T$, say $u_1$ and $u_2$, then $(u_1,x,u_2)$ is a path of
length 2 in $T$. By Lemma~\ref{lem5.1}, $b_{\rm R}(T)\leq
d_T(u_1)+d_T(x)+d_T(u_2)-3\leq \Delta(T)-1$.

Assume now that each vertex of $T$ is adjacent to at most one vertex
of degree one. Then $T$ has a vertex $u$ of degree 2 adjacent to
exactly one vertex, say $v$, of degree one. Let $w$ be the other
vertex adjacent to $u$. Then $(v,u,w)$ is a path of length 2 in $T$.
By Lemma~\ref{lem5.1}, $b_{\rm R}(T)\leq d_T(v)+d_T(u)+d_T(w)-3\leq
\Delta(T)$.
\end{pf}

\begin{lem}\label{lem5.2}
Let $G$ be a connected graph of order $n\ (\ge 3)$ and $\gamma_{\rm
R}(G)=\gamma(G)+1$. If there is an set $B$ of edges with
$\gamma_{\rm R}(G-B)=\gamma_{\rm R}(G)$, then
$\Delta(G)=\Delta(G-B)$.
\end{lem}

\begin{pf}
Since $G$ is connected and $n\geq 3$, $\gamma_{\rm
R}(G)=\gamma(G)+1\leq n-1$. Since $\gamma_{\rm R}(G-B)=\gamma_{\rm
R}(G)\leq n-1$,  $G-B$ is nonempty. By Lemma~\ref{lem2.5} and
Lemma~\ref{lem2.6}, $\gamma_{\rm R}(G-B)\geq \gamma(G-B)+1$. Since
$\gamma_{\rm R}(G-B)=\gamma_{\rm R}(G)=\gamma(G)+1\leq
\gamma(G-B)+1$, $\gamma_{\rm R}(G-B)=\gamma(G-B)+1$ and
$\gamma(G-B)=\gamma(G)$.

If $G-B$ is connected, then by Lemma~\ref{lem2.7},
$\Delta(G-B)=n-\gamma(G-B)=n-\gamma(G)=\Delta(G)$.

If $G-B$ is disconnected, then let $G_1$ be a nonempty connected
component of $G-B$. By Lemma~\ref{lem2.5} and~\ref{lem2.6},
$\gamma_{\rm R}(G_1)\geq \gamma(G_1)+1$. Then
$\gamma(G)+1=\gamma_{\rm R}(G-B)\ge \gamma_{\rm R}(G_1)+\gamma_{\rm
R}(G-G_1)\geq \gamma(G_1)+1+\gamma(G-G_1)\geq \gamma(G)+1$, thus
$\gamma_{\rm R}(G_1)=\gamma(G_1)+1$,  $\gamma_{\rm
R}(G-G_1)=\gamma(G-G_1)$ and $\gamma(G)=\gamma(G_1)+\gamma(G-G_1)$.
By Lemma~\ref{lem2.6}, $G-G_1$ is empty and hence
$\gamma(G-G_1)=|V(G-G_1)|$. By Lemma~\ref{lem2.7},
$\Delta(G_1)=|V(G_1)|-\gamma(G_1)=n-|V(G-G_1)|-\gamma(G_1)
=n-\gamma(G-G_1)-\gamma(G_1)=n-\gamma(G)=\Delta(G)$.
\end{pf}

\begin{thm}\label{thm5.3}
Let $G$ be a connected graph of order $n\ (\ge 3)$ and $\gamma_{\rm
R}(G)=\gamma(G)+1$. Then $b_{\rm R}(G)\leq \min
\{b(G),n_{\Delta}\}$, where $n_{\Delta}$ is the number of vertices
with maximum degree $\Delta$ in $G$.
\end{thm}

\begin{pf}
Since $n\geq 3$ and $G$ is connected, $\Delta(G)\geq 2$ and hence
$\gamma(G)\leq n-2$. Let $B$ be a minimum bondage set of $G$. Then
$G-B$ is nonempty and by Lemma~\ref{lem2.5} and Lemma~\ref{lem2.6}.
Thus, $\gamma_{\rm R}(G-B)\geq \gamma(G-B)+1>\gamma(G)+1=\gamma_{\rm
R}(G)$ and hence $b_{\rm R}(G)\le b(G)$.

We now prove that $b_{\rm R}(G)\leq n_{\Delta}$. By Lemma
\ref{lem2.7}, $\gamma_{\rm R}(G) =\gamma(G)+1$ if and only if there
is a vertex of degree $n-\gamma(G)$. If there is a vertex $s$ in $G$
such that $d_G(s)> n-\gamma(G)$, let $f(s)=2$ and $f(w)=1$ for any
vertex $w$ not in $N_G[s]$, then $f$ is a Roman dominating function
of $G$ with $f(G)=\gamma(G)$, a contradiction. Thus,
$\Delta(G)=n-\gamma(G)$. We can remove a smallest edge set $B$ with
$|B|\leq n_{\Delta}$ edges from $G$ such that
$\Delta(G-B)<\Delta(G)=n-\gamma(G)$ and $G-B$ is nonempty. Since
$G-B$ is nonempty,  by Lemma~\ref{lem2.5} and Lemma~\ref{lem2.6},
$\gamma_{\rm R}(G-B)\geq \gamma(G-B)+1$. Assume $\gamma_{\rm
R}(G-B)=\gamma_{\rm R}(G)$, then by Lemma~\ref{lem5.2},
$\Delta(G-B)=\Delta(G)=n-\gamma(G)$, a contradiction. Hence $b_{\rm
R}(G)\leq |B|\leq n_{\Delta}$.
\end{pf}

\begin{thm}\label{thm5.4}
For Roman graph $G$, $b_{\rm R}(G) \geq b(G)$.
\end{thm}

\begin{pf}
Let $B$ be a minimum Roman bondage set of $G$, then $\gamma_{\rm
R}(G-B)>\gamma_{\rm R}(G)=2\gamma(G)$. By Lemma~\ref{lem2.5},
$\gamma_{\rm R}(G-B)\leq 2\gamma(G-B)$, then $\gamma(G-B)>\gamma(G)$
and hence $b_{\rm R}(G) \geq b(G)$.
\end{pf}

The equality in Theorem~\ref{thm5.4} can hold, for example,
$b(C_{3k})=2=b_{\rm R}(C_{3k})$, and the strict inequality can also
hold, for example, $b(C_{3k+2})=2<3=b_{\rm R}(C_{3k+2})$.

\begin{thm}\label{thm5.5}
Let $G$ be a nonempty graph with $\gamma_{\rm R}(G)\geq 3$. Then
$b_{\rm R}(G)\leq (\gamma_{\rm R}(G)-2)\Delta(G)+1$.
\end{thm}

\begin{pf}
The proof proceeds by induction on $\gamma_{\rm R}(G)$.

We first assume that $\gamma_{\rm R}(G)=3$. Then by
Lemma~\ref{lem2.3}, $\Delta(G)=|V(G)|-2$. Assume that $b_{\rm
R}(G)\geq \Delta(G)+2$. Let $u$ be a vertex of maximum degree in
$G$. We have $\gamma_{\rm R}(G-u)=\gamma_{\rm R}(G)-1=2$. There is a
vertex $v$ that is adjacent to every vertex in $G-u$ and hence
$vu\notin E(G)$. Since $b_{\rm R}(G-u)\geq 2$, then for any edge
$e\in E_{G-u}(v)$, $\gamma_{\rm R}(G-u-e)=2$. Thus there is a vertex
$w$ that is adjacent to every vertex of $G-u-e$. But, since $v$ is
the only vertex of $G$ that is not adjacent to $u$, $wu\in E(G)$,
$d_G(w)=|V(G)|-1$, a contradiction. Thus, $b_{\rm R}(G)\leq
\Delta(G)+1$ if $\gamma_{\rm R}(G)=3$.

Assume the induction hypothesis for any integer $k$ and any graph
$H$ with $\gamma_{\rm R}(H)=k\geq 3$. Let $G$ be a nonempty graph
with $\gamma_{\rm R}(G)=k+1$, and assume that $b_{\rm R}(G)\geq
(k-1)\Delta(G)+2$. For any vertex $u$ of $G$, let $H=G-u$. Then,
$\gamma_{\rm R}(H)=\gamma_{\rm R}(G)-1=k$ since $d_G(u)<b_{\rm
R}(G)$. By the inductive hypothesis and by Lemma~\ref{lem5.1}, we
have
 $$
 \begin{array}{rl}
 b_{\rm R}(G)&\leq b_{\rm R}(H)+d_G(u)\\
 &\leq (k-2)\Delta(H)+1+d_G(u)\\
 &\leq (k-2)\Delta(G)+1+\Delta(G)\\
 &=(k-1)\Delta(G)+1,
 \end{array}
 $$
a contradiction. Thus, $b_{\rm R}(G)\leq  (k-1)\Delta(G)+1$, and by
the principle of mathematical induction, $b_{\rm R}(G)\leq
(\gamma_{\rm R}(G)-2)\Delta(G)+1$.
\end{pf}

Use $\kappa(G)$ (resp. $\lambda(G)$) to denote the
vertex-connectivity (resp. the edge-connectivity) of a connected
graph $G$ which is the minimum number of vertices (resp, edges)
whose removal results in $G$ disconnected. The famous Whitney's
inequality can be stated as
$\kappa(G)\leqslant\lambda(G)\leqslant\delta(G)$ for any graph $G$.
A subset $F\subseteq E(G)$ is called a $\lambda$-cut if
$|F|=\lambda(G)$ and $G-F$ is disconnected.

\begin{thm}\label{thm5.8}
If $G$ is a connected graph with order at least 3, then $b_{\rm
R}(G)\leq 2\Delta(G)+ \lambda(G)-3$, where $\lambda(G)$ is the
edge-connectivity of $G$.
\end{thm}

\begin{pf}
Let $G$ be a connected graph with edge-connectivity $\lambda(G)$ and
$F$ be $\lambda$-cut of $G$. Then $H=G-F$ has exact two connected
components. Let $x,y\in V(G)$, $xy\in F$, and $H_x$ and $H_y$ denote
the components of $G-F$ containing $x$ and $y$, respectively.
Without loss of generality, let $z$ be adjacent to $x$ in $H_x$
since $|V(G)|\geq 3$. Let $B=F\cup E_{H_x}(x)\cup E_{H_x}(z)-xz$ and
$f$ be a $\gamma_{\rm R}$-function of $G'=G-B$. Then $x$ and $z$ is
only adjacent to each other in $G'$, and so we can assume $f(x)=2$
and $f(z)=0$. We construct a Roman dominating function $f'$ of $G$
with $f'(G)<f(G')$.

If $V(H_y)=\{y\}$, then $f(y)=1$. Let $f'=f$ except $f'(y)=0$. Then
$f'$ is a Roman dominating function of $G$ with $f'(G)<f(G')$. Thus,
$b_{\rm R}(G)\leq |B|\leq 2\Delta(G)+ \lambda(G)-3$. In the
following, we assume $|V(H_y)|\geq 2$.

If $\gamma_{\rm R}(H_y-y)\geq \gamma_{\rm R}(H_y)$, then
 $$
 \gamma_{\rm R}(G-(F\cup E_{H_y}(y)))\geq \gamma_{\rm
 R}(H_x)+\gamma_{\rm R}(H_y)+1\geq \gamma_{\rm R}(G)+1.
 $$
Thus
 $$
 \begin{array}{rl}
 b_{\rm R}(G)&\leq |F\cup E_{H_y}(y))|\leq \Delta(G)+ \lambda(G)-1\\
 &\leq 2\Delta(G)+ \lambda(G)-3.
 \end{array}
 $$

If $\gamma_{\rm R}(H_y-y)= \gamma_{\rm R}(H_y)-1$, we can assume
that $f(y)=1$. Let $f'=f$ except $f'(y)=0$. Then $f'$ is a Roman
dominating function of $G$ with $f'(G)<f(G')$. Thus,
 $$
 b_{\rm R}(G)\leq |B|\leq 2\Delta(G)+ \lambda(G)-3.
 $$

The theorem follows.
\end{pf}

Considering vertex rather than edge-connectivity, we could
conjecture an analogy of Theorem ~\ref{thm5.8} by a similar
argument.

\begin{cjt}
If $G$ is a connected graph with order no less than 3, then $b_{\rm
R}(G)\leq 2\Delta(G)+ \kappa(G)-3$, where $\kappa(G)$ is the
vertex-connectivity of $G$.
\end{cjt}

\begin{thm}\label{thm5.9}
If $G$ is a nonempty graph with a unique minimum
Roman dominating function, then $b_{\rm R}(G)=1$.
\end{thm}

\begin{pf}
Let $f$ be the unique $\gamma_{\rm R}$-function on $G$, and let $x$
be a vertex in $G$ with $f(x)=0$. Then there is a vertex $y\in
N_G(x)$ with $f(y)=2$. If there are at least two vertices $y,z\in
N_G(x)$ such that $f(y)=f(z)=2$ for each vertex $x$ with $f(x)=0$.
Then let $f'=f$ except that $f'(x)=2$ and $f'(y)=0$ and $f'$ is a
$\gamma_{\rm R}$-function on $G$ as well, which is a contradiction
to the uniqueness of $f$. Thus, there is a unique $y\in N_G(x)$ with
$f(y)=2$ for a vertex $x$ with $f(x)=0$. Then $\gamma_{\rm
R}(G-xy)>\gamma_{\rm R}(G)$, which implies that $b_{\rm R}(G)=1$.
\end{pf}

\begin{thm}\label{thm5.10}
If $G$ is a vrc-graph with $\gamma_{\rm R}(G)=3$, then $b_{\rm
R}(G)\leq \Delta(G)+1$.
\end{thm}

\begin{pf}
By Lemma~\ref{lem2.9}, $G$ is a vc-graph with $\gamma(G)=2$. By
Lemma~\ref{lem2.8}, $G$ is a complete $K_{2t} (t\geq 2)$ with a
perfect matching $M$ removed. Thus, $G$ is $\Delta(G)$-regular,
where $\Delta(G)=2t-2$. Let $uv\in M$. Then $v$ is the only vertex
not adjacent to $u$ in $G$. Let $H=G-u$. Then $\gamma_{\rm R}(H)=2$
since $G$ is a vrc-graph with $\gamma_{\rm R}(G)=3$. Note that the
vertex $v$ is the only vertex adjacent to all the other vertices in
$H$ adjacent to each of other vertices in $H$. Thus $H$ has a unique
minimum Roman dominating function $f$ with $f(v)=2=\gamma_{\rm
R}(H)$. By Theorem~\ref{thm5.9}, $b_{\rm R}(H)=1$ and hence $b_{\rm
R}(G)\leq \Delta(G)+1$.
\end{pf}

\begin{thm}
If there exists at least one vertex $u$ in a graph $G$ with
$\gamma_{\rm R}(G-u)\geq \gamma_{\rm R}(G)$, then $b_{\rm
R}(G)=d_G(x)\leq \Delta(G)$.
\end{thm}

\begin{pf}
Since $\gamma_{\rm R}(G-E_G(u))=\gamma_{\rm R}(G-u)+1>\gamma_{\rm
R}(G)$, $b_{\rm R}(G)=d_G(x)\leq \Delta(G)$.
\end{pf}

\begin{cor}
Let $G$ be a graph of order $n$. If $\gamma_{\rm R}(G)=3\neq n$,
then $b_{\rm R}(G)\leq \Delta+1$.
\end{cor}

\begin{prob}
Whether or not there exits a positive integer $c$ such that $b_{\rm
R}(G)\leq \Delta(G)+c$ for any graph $G$ of order $n$ and
$\gamma_{\rm R}(G)\neq n$.
\end{prob}

The {\it vertex covering number} $\beta(G)$ of $G$ is the minimum
number of vertices that are incident with all edges in $G$. If $G$
has no isolated vertices, then $\gamma_{\rm R}(G)\leq 2\gamma(G)
\leq 2\beta(G)$. If $\gamma_{\rm R}(G)=2\beta(G)$, then $\gamma_{\rm
R}(G)=2\gamma(G)$ and hence $G$ is a Roman graph. In~\cite{v94},
Volkmann gave a lot of graphs with $\gamma(G)=\beta(G)$.

\begin{thm}\label{thm5.11}
Let $G$ be a graph with $\gamma_{\rm R}(G)=2\beta(G)$. Then\\
(1) $b_{\rm R}(G)\geq \delta(G)$;\\
(2) $b_{\rm R}(G)\geq \delta(G)+1$ if $G$ is a vrc-graph.
\end{thm}

\begin{pf}
Let $G$ be a graph with $\gamma_{\rm R}(G)=2\beta(G)$.

(1) Without loss of generality, Assume $\delta(G)\geq 2$. Let
$B\subseteq E(G)$ with $|B|\leq \delta(G)-1$. Then $\delta(G-B)\geq
1$ and so $\gamma_{\rm R}(G)\leq \gamma_{\rm R}(G-B) \leq
2\beta(G-B) \leq 2 \beta(G)=\gamma_{\rm R}(G)$. Thus, $B$ is not a
Roman bondage set of $G$, and so $b_{\rm R}(G)\geq \delta(G)$.

(2) From the above proof, every Roman bondage set $B$ contains at
least all edges incident with some vertex $x$, so that $G-B$ has an
isolated vertex. On the other hand, if $G$ is a vrc-graph, then
$\gamma_{\rm R}(G-x) <\gamma_{\rm R}(G)$ for any vertex $x$, which
implies that the removal of all edges incident with $x$ can not
enlarge the Roman domination number. Hence $b_{\rm R}(G)\geq
\delta(G)+1$.
\end{pf}

\end{document}